\documentclass{elsarticle}
\usepackage{graphicx} % Required for inserting images
\usepackage{amsmath}
\usepackage{amssymb}
\usepackage{xcolor}
\usepackage{subcaption}

\journal{Journal of Computational Physics}

\begin{document}

\begin{frontmatter}

\title{Numerical Robustness of PINNs for Multiscale Transport Equations}

\begin{abstract}
We investigate the numerical solution of multiscale transport equations using Physics Informed Neural Networks (PINNs) with ReLU activation functions. Therefore, we study the analogy between PINNs and Least-Squares Finite Elements (LSFE) which lies in the shared approach to reformulate the PDE solution as a minimization of a quadratic functional. We prove that in the diffusive regime, the correct limit is not reached, in agreement with known results for first-order LSFE. A diffusive scaling is introduced that can be applied to overcome this, again in full agreement with theoretical results for LSFE. We provide numerical results in the case of slab geometry that support our theoretical findings.
\end{abstract}
%%Research highlights

\author{Alexander Jesser\corref{cor1}\fnref{affiliation1,affiliation2}}
\ead{alexander.jesser@kit.edu}
\author{Kai Krycki\fnref{affiliation3}}
\ead{krycki@fh-aachen.de}
\author{Ryan G. McClarren\fnref{affiliation4}}
\ead{rmcclarr@nd.edu}
\author{Martin Frank\fnref{affiliation1}}
\ead{martin.frank@kit.edu}
\cortext[cor1]{Corresponding author}

\affiliation[affiliation1]{organization={Karlsruhe Institute of Technology, Scientific Computing Center (SCC)},
            addressline={Hermann-von-Helmholtz-Platz 1},
            postcode={76344},            
            city={Eggenstein-Leopoldshafen},
            country={Germany}}

\affiliation[affiliation2]{organization={Aachen Institute for Nuclear Training GmbH},
            addressline={Cockerillstrasse 100 (DLZ)},
            postcode={52222},            
            city={Stolberg(Rhld.)},
            country={Germany}}

\affiliation[affiliation3]{organization={FH Aachen - University of Applied Sciences, Department of Aerospace Engineering},
            addressline={Hohenstaufenallee 6},
            postcode={52064},            
            city={Aachen},
            country={Germany}}

\affiliation[affiliation4]{organization={University of Notre Dame, Department of Aerospace and Mechanical Engineering},
            addressline={257 Fitzpatrick Hall},
            postcode={46556},            
            city={Notre Dame},
            country={USA}}            

\end{frontmatter}

\section{Introduction}\label{sec:intro}
Physics Informed Neural Networks (PINNs) have recently attracted a great deal of attention. First proposed in \cite{PINNs-org} in 2019, PINNs provide a method to solve partial differential equations (PDEs) using (deep) neural networks (NN) by incorporating the corresponding equation in a loss function in a training process. Successively, points in the computational domain of the PDE are chosen and the weights of the NN are adapted to solve the PDE pointwise, resulting in an approximation of the solution on the complete computational domain. In that way, the numerical solution of the PDE corresponds to training of a NN. First computational packages are available which therefore provide easy to use black-box solver for PDEs, what makes the method attractive for users.\\

Although already heavily used in various applications \cite{PINNs-nature,Cai2021,Cai2021-FM,Cuomo2022,Mishra2021}, only a limited number of results regarding the numerical analysis of PINNs exist so far \cite{DeRyck_Mishra_2024}. For the classes of linear second order elliptic and parabolic PDEs convergence of the discretized solution to the PDE solution can be shown in a strong sense \cite{PINN_convergence2020}. A more detailed numerical analysis of PINNs (convergence, error bounds, stability etc.) needs to be developed based on these first results in order to gain a deeper understanding of this undoubtedly powerful and useful method.\\
With this paper, we contribute to the numerical analysis of PINNs by studying an analogy between PINNs and the Least Squares Finite Element method (LSFE). For a specific class of PDEs (neutron transport equations) we study similarities based on an existing theory for LSFE for numerical solution in the so-called diffusive regime (cf. \cite{Manteuffel1993}, \cite{Manteuffel1998}, \cite{Manteuffel2000}). The LSFE is a numerical approach used to solve partial differential equations and involves reformulating the problem into a minimization of a least-squares functional. A numerical solution is obtained by solving for a finite element representation of the corresponding variational formulation.\\

The analogy between PINNs and the LSFE lies in their shared approach of reformulating the solution of differential equations as the minimization of a (typically quadratic) functional. In LSFE, this involves the construction of a least-squares functional based on the residual of the differential equation, which, when minimized, yields the solution to the problem in a variational form. Similarly, PINNs define a loss function based on the residuals of the differential equation and corresponding boundary conditions. Minimizing this loss function drives the neural network to approximate the solution. Hence, although not explicitly stated in variational terms, the minimization of the loss function in PINNs can be viewed as finding an approximate solution that satisfies the differential equation in an integral sense across the domain. Using an activation function such as ReLU for PINNs, the corresponding neural network can be interpreted as a piecewise linear approximation of the solution comparable to the finite element representation.\\

For this reason, we expect both methods to show comparable behavior in situations where numerically resolving multiscale effects in a PDE becomes difficult. A prominent and widely studied example is neutron transport in the so-called diffusive regime, which is briefly introduced in the following.

\subsection*{Neutron Transport in the Diffusive Regime}
The transport of neutrons in a general background medium is governed by the linear Boltzmann Equation \cite{cerc1988}. For a time-independent source-term in the case of mono-energetic neutrons and isotropic scattering the steady-state form reads
\begin{align}
    \mathbf{\Omega} \cdot \nabla \psi(\mathbf{x}, \mathbf{\Omega}) + \sigma_t \psi(\mathbf{x}, \mathbf{\Omega}) = \frac{\sigma_s}{4\pi} \int_{4\pi} \psi(\mathbf{x}, \mathbf{\Omega}') \, d\mathbf{\Omega}' + Q(\mathbf{x}, \mathbf{\Omega}),
\end{align}
where $\psi(\mathbf{x}, \mathbf{\Omega})$ is the angular flux of neutrons at position $\mathbf{x}$ and direction $\mathbf{\Omega}$,     $\mathbf{\Omega} \cdot \nabla$ is the transport operator, $\sigma_t$ denotes the total and $\sigma_s$ the scattering macroscopic cross-section. $S(\mathbf{x}, \mathbf{\Omega})$ is the source term at position $\mathbf{x}$ in direction $\mathbf{\Omega}$ and the integral $\frac{\sigma_s}{4\pi} \int_{4\pi} \psi(\mathbf{x}, \mathbf{\Omega}') \, d\mathbf{\Omega}'$ represents the isotropic scattering contribution. \\

\begin{figure}[htbp]
\centering
  \centering
\includegraphics[width=0.75\textwidth]{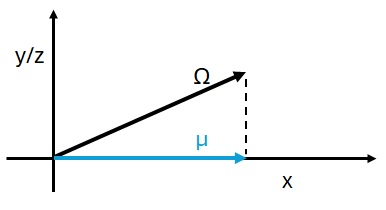}
\caption{In \textit{slab geometry}, transport is projected onto the x-axis. The problem is symmetric in the y/z-plane.}\label{fig:slab_geometry}
\end{figure}

To simplify the analysis, it is common to investigate the \textit{slab geometry} case \cite{Lewis_Miller_1984} (one space and directional dimension), where all relevant effects in the diffusive regime can be shown. Here, the medium is assumed to be infinite in two dimensions (e.g., the $y$- and $z$-directions) and finite in the $x$-direction. The medium is therefore bounded by two parallel planes which form a slab through which the neutrons travel. In the following, we denote the $x$-coordinates of the planes that limit the slab by $x_l$ for the left plane and $x_r$ for the right plane. The (one-dimensional) spatial coordinate $x$ is used to describe the position of a neutron within the slab, the direction of neutron movement becomes a one-dimensional angle $\mu=\cos(\theta)$, where $\theta$ is the angle between the neutron’s direction of travel and the $x$-axis. Fig. \ref{fig:slab_geometry} depicts the \textit{slab geometry} case.\\

In the following, we assume an isotropic source, i.e. $Q$ independent of $\mu$. The steady state, single energy group neutron transport equation in slab geometry is then formulated as
\begin{align}\label{eq:nte}
    L\Psi(x,\mu)=Q(x)
\end{align}
with source term $Q$ and operator
\begin{align}\label{eq:transport_operator}
    L\Psi(x,\mu)=\left(\mu\frac{\partial}{\partial x}+\sigma_t(I-P)+\sigma_a P\right)\Psi(x,\mu),
\end{align}
where $I$ denotes the identity operator and $P$ the projection onto the space of all $L^2$-functions independent of angle
\begin{align}
    (P\psi)(x)=\frac{1}{2}\int_{-1}^{1} \psi(x,\mu) d\mu.
\end{align}
Boundary conditions can be introduced by a boundary operator $B$, so that
\begin{align}\label{eq:boundary_operator}
    B\Psi(x,\mu)=g(x,\mu)~\forall (x,\mu)\in \Gamma^-,
\end{align}
with $\Gamma^-=\{(x,\mu)\in\partial D\times [-1,1]:n(x)\cdot\mu< 0\}$, where $\partial D$ is the boundary of the computational domain $D$ and $n$ is the outer normal vector. The boundary conditions in this form specify the flux distribution entering $D$ over $\partial D$. We introduce the boundary conditions in detail later in Section \ref{sec:num-results}.\\
The multiscale effects in this setting originate from two different scales relevant to this problem: the mean-free path is small compared to the size of the physical domain. The numerical instability of PINNs that we discuss in this paper is a consequence of the different scales and the stiffness they induce. In the diffusive regime scattering dominates the effects of free transport or absorption and the mean-free path of particles is small. The diffusive regime can typically be characterized by a small parameter $\varepsilon>0$, which can be interpreted as the ratio of the mean-free path to the physical size of the domain of computation. Introducing such $\varepsilon>0$ we scale the cross sections and source term accordingly \cite{larsen1992} as 
\begin{align}
    \sigma_t\rightarrow\frac{1}{\varepsilon},~\sigma_a\rightarrow \alpha \varepsilon,~Q\rightarrow\varepsilon Q,
\end{align}
where $\alpha$ is assumed to be $\mathcal{O}(1)$. By applying this scaling on Eq. (\ref{eq:nte}) we obtain 
\begin{align}\label{eq:nte_diff_limit}
    L_\varepsilon\Psi^\varepsilon(x,\mu)=\varepsilon Q(x)
\end{align}
where the neutron flux in the diffusion limit is denoted as $\Psi^\varepsilon$ and the scaled operator $L_\varepsilon$ reads
\begin{align}\label{L_diff_limit}
    L_\varepsilon\Psi^\varepsilon(x,\mu)=\left(\mu\frac{\partial}{\partial x}+\frac{1}{\varepsilon}(I-P)+\varepsilon\alpha P\right)\Psi^\varepsilon(x,\mu).
\end{align}
It can be shown \cite{larsen1980} that in the limit of $\varepsilon\rightarrow 0$ the solution of the transport equation $\Psi^\varepsilon$ converges to a function $\phi_0$, which is independent of $\mu$ and the solution of an associated diffusion equation
\begin{align}
    -\frac{\partial}{\partial x}\frac{1}{3\sigma_t}\frac{\partial \phi_0}{\partial x}(x)+\sigma_a\phi_0(x)=Q(x).
\end{align}
For numerical methods it is therefore of key importance to reproduce this behavior on a discretized level. A method that ensures this is called asymptotic preserving \cite{Hu2017}. For the LSFE method a corresponding analysis has been carried out (cf. \cite{Manteuffel1993}, \cite{Manteuffel1998}, \cite{Manteuffel2000}) and we build upon these results for the investigation of PINNs. A related analysis using PINNs for multiscale time dependent linear transport equations has been carried out in \cite{ShiJin_2023}. Using LSFE without additional scaling (see Section \ref{sec:main-theorem}), the numerical solution produces an incorrect limit as $\varepsilon\to 0$. We show in this paper that the same happens when PINNs are applied to solve this class of equations and that the same scaling as derived for LSFE can be applied to ensure a correct limiting behavior for PINNs, underlining the close connection between the two numerical methods. \\

The paper is structured as follows. Sections \ref{sec:PINNs} and \ref{sec:LSFE} introduce PINNs and LSFE respectively as numerical methods for the solution of PDEs. Here, we also introduce the notations of the methods that we use in the following. The main result is stated in Section \ref{sec:main-theorem}. We prove that numerical solutions of PINNs in the diffusive regime do not converge to the correct solution, and show that a diffusive scaling can be applied to overcome this, in full agreement with the theoretical results for LSFE. Numerical results are presented in Section \ref{sec:num-results}, and consequences of our results as well as future research directions are discussed in the final Section \ref{sec:discussion}.

%-------------- SECTION ---------------------------------------------------------------------------------------
\section{Physics Informed Neural Networks}\label{sec:PINNs}

Physics Informed Neural Networks, PINNs for short, are a numerical method for solving general partial differential equations \cite{PINNs-org,PINNs-nature}. A neural network serves as an interpolation function for the solution of the PDE. It is obtained by successively selecting points in the computational domain and adjusting the parameters of the neural network to obtain the desired output, which is commonly referred to as training the NN. The method is therefore mesh-free and can serve as a black-box solver for a PDE, since only the form of the equation needs to be known. The mathematical formulation of the method is briefly presented below.\\
In the following, we use a notation similar to \cite{PINN_convergence2020}. Let $\Psi^M(x,\mu):\mathbb{R}^{2} \rightarrow \mathbb{R}^{1}$ be a feed-forward neural network (FNN) with $M$ layers. Then its output is recursively defined by
\begin{align}\label{eq:NN_recursice_definition}
    &\Psi^0(x,\mu)~=(x,\mu)\in\mathbb{R}^{2}\\ \nonumber
    &\Psi^m(x,\mu)=W^m\sigma(\Psi^{m-1}(x,\mu))+b^m~~\in\mathbb{R}^{M_l},~2\leq m \leq M-1\\ \nonumber
    &\Psi^M(x,\mu)=W^M\sigma(\Psi^{M-1}(x,\mu))+b^M~\in \mathbb{R}^{1} 
\end{align}
with the weight matrices $W^m\in\mathbb{R}^{N_m,N_{m-1}}$, the biases $b^m\in\mathbb{R}^{N_m}$ and the activation function $\sigma$. Let $\Vec{n}=\{N_0,...,N_M\}\in \mathbb{N}^M$ with $\theta=\{W^m,b^m\}_{1 \leq m\leq M}$ be the parameter space of the network. Since a neural network $\Psi^M(x,\mu)$ depends on the architecture $\Vec{n}$ and its parameters $\theta$, we can denote it as $\Psi^M=\Psi^M(x,\mu,\Vec{n},\theta)$. For a fixed $\Vec{n}$, we write $\Psi^M(x,\mu,\theta)$.\\
A numerical solution of a PDE using such a FNN can be obtained by incorporating the PDE residual into the loss function during the training process, thus obtaining a 'physics informed' neural network. In the following, we assume that a unique classical solution of the PDE exists. The PDE and its boundary conditions are given by
\begin{align}
    L[\Psi](x,\mu)=Q(x)~\forall x\in D,~B[\Psi](x,\mu)=g(x,\mu)~\forall (x,\mu)\in \Gamma^-.
\end{align}
Here $L$ denotes a differential operator, $Q$ is a source term and $B$ is a boundary operator. Hence, the PINN maps the spatial variable $x$ to an approximate solution of the PDE via $\Psi^M(x,\mu)\approx \Psi(x,\mu)$. The loss function of the neural network is given by
\begin{align}\label{eq:loss_function}
    \text{Loss}[\Psi]=\frac{1}{|T_D|}\sum_{(x,\mu)\in T_D}(L[\Psi](x,\mu)-Q(x))^2,
\end{align}
where $T_D$ denotes a set of training points within the domain $D$. This set is often generated using a quasi-random low-discrepancy sequence like the Sobol sequence \cite{Sobol1967}. The training process of the PINN is the equivalent to solving the PDE in classical methods. It is a minimization of the loss function in the set of functions $V$ which is defined by
\begin{align}
    V = \{\Psi^M(\cdot,\cdot,\theta):\mathbb{R}^{2} \rightarrow \mathbb{R}^{1}|\theta=\{W^m,b^m\}_{1 \leq m\leq M} \}.
\end{align}
In general, $V$ is not a vector space. The boundary conditions are incorporated by a restriction of the set of functions so that
\begin{align}
    V_b = \{\Psi^M(\cdot,\cdot,\theta):\mathbb{R}^{2} \rightarrow \mathbb{R}^{1}|\theta=\{W^m,b^m\}_{1 \leq m\leq M}\nonumber\\~\text{where}~B\Psi^M(x,\mu)=g(x,\mu)~\forall (x,\mu)\in \Gamma^- \}.
\end{align}
Alternatively, the boundary condition can be included in the loss function. We will describe this approach in Subsec. \ref{subsec:implementation}. For now, we can write the training process of the PINN as a minimization of the loss function in the set of functions $V_b$:
\begin{align}
    \Psi_{min} = \min_{\Psi\in V_b}\text{Loss}[\Psi].
\end{align}
Equivalently, we can interpret it as a minimization of the loss function in the parameter space with gradient-based optimizers such as Adams \cite{Kingma_Ba_2015} or L-BFGS \cite{L-BFGS}:
\begin{align}\label{eq:PINN_nonlinear_problem}
    \Psi_{min} = \min_{\theta}\text{Loss}[\Psi](\theta).
\end{align}
It is important to note that $\Psi^M(x,\mu)$ using the rectified linear unit (ReLU) as activation function is piecewise linear. Moreover, it can be shown that Neural Networks with ReLU activation functions and sufficiently many layers can reproduce all piecewise linear basis functions from the finite element method \cite{RELU_DNN_and_LFE}.\\

The problem considered in this work depends, in principle, on the architecture of the underlying NN, including the number of layers and neurons. However, our analysis remains valid regardless of the specific architecture, provided that a feed-forward NN is used, which is sufficiently large to accurately approximate the solution to the PDE. The theoretical analysis focuses on networks with ReLU activation functions, while the case of hyperbolic tangent activation functions is also explored in the numerical experiments. In general, using different activation functions (or a combination of them) or Non-Feedforward Neural Networks might lead to different results in the diffusion limit.

%-------------- SECTION ---------------------------------------------------------------------------------------
\section{Least-Squares Finite Elements}\label{sec:LSFE}
The Least-Squares Finite Element (LSFE) method relies on a finite element discretization of a variational formulation of the underlying PDE. The least-squares formulation of Eq. (\ref{eq:nte}) (neutron transport in slab geometry) is given by
\begin{align}
     \min_{\Psi\in U} F[\Psi]~~\text{with}~~F[\Psi]=\int_{x_l}^{x_r}\int_{-1}^{1} (L\Psi(x,\mu)-Q(x,\mu))^2 d\mu dx
\end{align}
where the boundary condition (Eq. \ref{eq:boundary_operator}) is incorporated by a restriction of the function space:
\begin{align}
    U_b=\{u(x,\mu)\in L^2(D):\mu \frac{\partial u}{\partial x}\in L^2(D)~\text{where}~Bu(x,\mu)=0~\forall (x,\mu)\in \Gamma^- \}
\end{align}
For $\Psi$ to be a minimizer of $F$, it is a necessary condition that the first variation vanishes for all admissible $u\in U_b$, resulting in the problem: find $\Psi\in U$ such that
\begin{align}\label{eq: a_psi_l_formal}
    a(\Psi,u)=l(u)
\end{align}
for all admissible $u\in U_b$ with the bilinear form
\begin{align}\label{eq:LSFE_var_form}
    a(\Psi,u):=\int_{x_l}^{x_r}\int_{-1}^{1} L\Psi(x,\mu)Lu(x,\mu) d\mu dx 
\end{align}
and the functional
\begin{align}
    l(u)=\int_{x_l}^{x_r}\int_{-1}^{1} Q(x,\mu)Lu(x,\mu) d\mu dx.
\end{align}
The system is then discretized by replacing the function space $U$ with a finite dimensional subspace $U_h\subset U$. Consequently, we obtain a finite dimensional linear system for the discretized neutron flux $\Psi_h$:
\begin{align}\label{eq:LSFE_discrete}
    a(\Psi_h,u_h)=l(u_h)~ \forall u_h\in U_h.
\end{align}
If $\{\lambda_p\}$ is a basis of $U_h$, we can represent the solution as an expansion in said basis:
\begin{align}
    \Psi_h=\sum_p c_p\lambda_p,
\end{align}
where $c_p$ are the expansion coefficients. This allows us to write Eq. (\ref{eq:LSFE_discrete}) as a linear system
\begin{align}\label{eq:LSFE_linear_problem}
    A\Psi_\lambda=b
\end{align}
where the components of $A$ and $l$ are given by $A=a(\lambda_p,\lambda_q)_{pq}$ and $b=l(\lambda_p)_p$, while $\Psi_\lambda$ is the vector of expansion coefficients $c_p$.\\
For the spatial discretization in the variable $x$ first-order Lagrange polynomials are used as basis functions. The discretization in the angular variable $\mu$ is with the method of moments is described in Subsec. \ref{subsec:angular_discretization}.\\
A general introduction to the least-squares finite element methods can be found in \cite{Bochev_Gunzburger}. Details about least-squares finite element solutions of the neutron transport equation in the diffusive regime are discussed in \cite{Manteuffel1998}.\\
A key difference between LSFE and PINNs is that, after parametrization, LSFE requires the solution of a linear problem (Eq. \ref{eq:LSFE_linear_problem}), whereas PINNs requires the solution a nonlinear optimization problem (Eq. \ref{eq:PINN_nonlinear_problem}). This means that for the PINN we look for weights and biases, and for LSFE we look for basis expansion coefficients.

%-------------- SECTION ---------------------------------------------------------------------------------------
\section{PINNs in the asymptotic limit}\label{sec:main-theorem}
\begin{figure}[htbp]
\centering
  \centering
\includegraphics[width=1.0\textwidth]{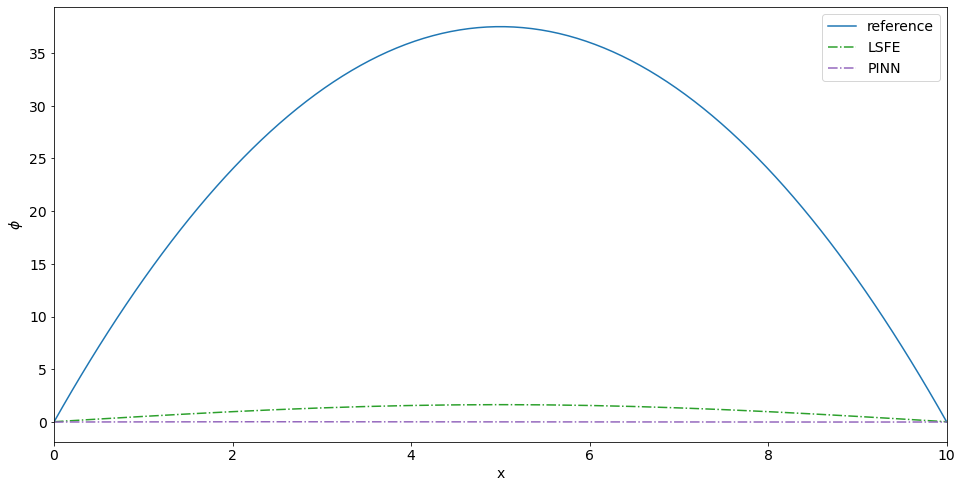}
\caption{Results for PINN with ReLU activation functions and first order LSFE in the asymptotic limit ($\varepsilon=10^{-4}$) compared to a MC reference solution. The setting is discussed in detail in Subsec. \ref{subsec:diff_limit}.}\label{fig:erronous_testcase}
\end{figure}

In the diffusive limit $\varepsilon\to 0$ as introduced in Section \ref{sec:intro} the neutron transport equation becomes singular and numerical methods often fail to produce valid approximations in this regime. Fig. \ref{fig:erronous_testcase} shows the behavior of a numerical solution obtained with a PINN with ReLU activation functions and a first-order LSFE solution in the asymptotic limit ($\varepsilon=10^{-4}$), compared to a Monte-Carlo (MC) reference solution. This testcase and its setting will be discussed in detail in Subsec. \ref{subsec:diff_limit}.

It can be seen that the neutron fluxes resulting from both numerical methods are close to zero and far away from the reference solution. The numerical analysis of the LSFE method yields an explanation to why this behavior occurs, and the methodological similarities between LSFE and PINNs suggested a corresponding behavior for PINNs.\\
In the following, we provide a numerical analysis that shows that a PINN with ReLU activation functions does in general not yield a correct solution in the diffusion limit.\\

We can express the solution of the transport equation by an expansion in the Legendre moments in angle. For $n\in\mathbb{N}$, let $P_n(\mu)$ be the normed n-th order Legendre polynomial on the interval $[-1,1]$. In this section, we use the normed Legendre polynomials defined as $p_n(\mu)=\sqrt{2n+1}P_n(\mu)$ to simplify the coefficients. Let
\begin{align}\label{eq: Legendre_moments}
    \phi_n(x) = \int_{-1}^1 p_n(\mu)\psi(x,\mu) d\mu
\end{align}
be the $n$-th angular moment of $\psi$. Then the expansion of the neutron flux in angle is given by
\begin{align}\label{eq: Legendre_expansion}
    \Psi(x,\mu)=\sum_{i=0}^{\infty}p_n(\mu)\phi_n(x).
\end{align}
Furthermore, we note that in the diffusion limit, it is sufficient to take only the first two  moments of the expansion into account \cite{Case_Zweifel_1967}.\\
In this situation we consider the set of functions
\begin{align}\label{eq:v_h}
    V_h=\{v_h\in V_b: v_h(x,\mu)=\phi_0(x)+\mu \phi_1(x), \text{where}~\phi_0,\phi_1\in \mathcal{P}(T_h)\},
\end{align}
where $\mathcal{P}(T_h)$ denotes the space of piecewise linear functions on the partition $T_h$ of the slab. Here the term partition corresponds to an interval within the slab, ie.e $[x_1,x_2]\subset[x_l,x_r]$. In the diffusive regime, the training of a PINN corresponds to the solution of the following minimization problem:
\begin{align}\label{eq:PDE_res_discrete}
   \Psi^\varepsilon_{min} = \min_{\Psi^\varepsilon\in V_h} \frac{1}{|T_D|} \sum_{i,j} (L_{\varepsilon}\Psi^\varepsilon(x_i,\mu_j)-\varepsilon Q(x_i,\mu_j))^2
\end{align}
In the following, we use the fact that for all $\xi>0$ it is possible to find a number of training points $|T_D(\xi)|$ so that the absolute difference between the discrete sum divided by the number of training points and the integral is smaller than the parameter $\xi$: 
\begin{align}\label{eq:PDE_res_continous}
   \bigg|\frac{1}{|T_D(\xi)|} \sum_{i,j} (L_{\varepsilon}\Psi^\varepsilon(x_i,\mu_j)-\varepsilon Q(x_i,\mu_j))^2 - \int_{x_l}^{x_r}\int_{-1}^{1} (L_\varepsilon\Psi^\varepsilon(x,\mu)-\varepsilon Q(x,\mu))^2 d\mu dx\bigg| < \xi.
\end{align}
This is ensured by the sequence we use to construct the training points. The training points are generated with a Sobol sequence which is constructed in a way that the sum divided by the number of (training) points converges to the integral \cite{Sobol1967}. Hence, if we replace in the following the term 'minimum' by '$\xi$-suboptimal solution', e.g. the $L^2$-norms difference of the minimum and the $\xi$-suboptimal solution is smaller than $\xi$, the results remain valid. For the sake of readability, in the following, we neglect this sublety and we assume that a function $\Psi^\varepsilon_{min}$ that minimizes Eq. (\ref{eq:min_F}) also minimizes Eq. (\ref{eq:PDE_res_discrete}) and use the term 'minimum'. Hence, instead of Eq. (\ref{eq:PDE_res_discrete}) we directly consider the following continuous minimization problem:
\begin{align}\label{eq:min_F}
     \min_{\Psi\in V_h} F[\Psi]~~\text{with}~~F[\Psi]=\int_{x_l}^{x_r}\int_{-1}^{1} (L_\varepsilon\Psi^\varepsilon(x,\mu)-\varepsilon Q(x,\mu))^2 d\mu dx.
\end{align}
To summarize, in the following analysis $\Psi^\varepsilon_{min}$ denotes the neural network trained on a Sobol sequence $T_D$ being a solution of Eq. (\ref{eq:PDE_res_discrete}). This neural network is an approximate solution of the neutron transport equation and is furthermore a piecewise linear, hence integrable function. By making use of the property of $\xi$-suboptimality we also consider $\Psi^\varepsilon_{min}$ to minimize the continuous functional Eq. (\ref{eq:min_F}). 

\subsection{Characterization of the minimizer}\label{subsec:minimizer}
In the following, we formally expand the solution $\Psi$ in powers of $\varepsilon$ and derive relations between the expansion coefficients. We then use these relations to prove that the PINN method with ReLU activation functions yields an incorrect solution in the diffusion limit. For our analysis we build on the theory developed in \cite{Ressel_phd_thesis} for LSFE, which we adapt for PINNs.\\

\newtheorem{thm}{Theorem}
\newtheorem{lem}{Lemma}
\begin{lem}
Let the Functional $F$ and the set of functions $V_h$ be as given in Eqs. (\ref{eq:min_F}) and (\ref{eq:v_h}). Suppose $\Psi^\varepsilon_{min}$ minimizes $F$ restricted to $V_h$. Suppose further that $\varepsilon\leq 1$ and that $\Psi^\varepsilon_{min}$ has an expansion
\begin{align}\label{eq:psi_min_expansion}
    \Psi^\varepsilon_{min}(x,\mu)=\phi_0^\varepsilon(x)+\mu \phi_1^\varepsilon(x)
\end{align}
with
\begin{align}\label{phi_developed}
    \phi_0^\varepsilon(x)=\sum_{\nu=0}^{\infty}\varepsilon^\nu \eta_\nu(x)~~~\text{and}~~~\phi_1^\varepsilon(x)=\sum_{\nu=0}^{\infty}\varepsilon^\nu \delta_\nu(x),
\end{align}
where $\eta_\nu$ and $\delta_\nu$ are independent of $\epsilon$. Then we have:
\begin{align}\label{main_theorem_relations}
    \delta_0(x)=0~~~\text{and}~~~\eta_0'(x) = -\delta_1(x).
\end{align}
\end{lem}

\newproof{pf}{Proof}
\begin{pf}
By inserting Eq. (\ref{eq:psi_min_expansion}) into Eq. (\ref{L_diff_limit}) and using that $P$ is a projection onto the zeroth moment ($P\Psi=\phi_0$), while $I-P$ projects onto moment one and all higher order moments, we obtain
\begin{align}\label{L_developed_psi}
    L_\varepsilon\Psi^\varepsilon_{min} = \mu\frac{\partial\phi_0^\varepsilon}{\partial x} + \mu^2 \frac{\partial \phi_1^\varepsilon}{\partial x} + \frac{1}{\varepsilon}\mu\phi_1^\varepsilon + \varepsilon\alpha\phi_0^\varepsilon
\end{align}
Inserting Eq. (\ref{phi_developed}) and sorting the terms in the order of $\varepsilon$ yields
\begin{align}\label{L_developed_psi_developed}
    L_\varepsilon\Psi^\varepsilon_{min} = \frac{1}{\varepsilon}[\mu\delta_0]+[\mu \eta_0' + \mu^2 \delta_0'+\mu\delta_1]+\mathcal{O}(\varepsilon).
\end{align}
By inserting Eq. (\ref{L_developed_psi_developed}) into Eq. (\ref{eq:min_F}), we obtain a representation of $F$ as a power series
\begin{align}
    F(\Psi^\varepsilon_{min})=\sum_{\nu=-2}^\infty \varepsilon^\nu F_\nu(\Psi_h)
\end{align}
with
\begin{align}
    &F_{-2}(\Psi^\varepsilon_{min}) = \int_{x_l}^{x_r} \int_{-1}^1 \mu^2\delta_0^2(x) d\mu dx\nonumber \\
    &~~~~~~~~~~~~~~=\frac{2}{3}\int_{x_l}^{x_r} \delta_0^2(x)dx
\end{align}
and
\begin{align}
    &F_{-1}(\Psi^\varepsilon_{min}) = 2\int_{x_l}^{x_r} \int_{-1}^1 \mu^2\eta_0'(x)\delta_0(x) + \mu^3\delta_0\delta_0' +\mu^2\delta_0(x)\delta_1(x) d\mu dx\nonumber\\
    &~~~~~~~~~~~~~~=\frac{4}{3}\int_{x_l}^{x_r} \eta_0'(x)\delta_0(x) + \delta_0(x)\delta_1(x) dx.
\end{align}
Since $0\in V_h$ and $\Psi^\varepsilon_{min}$ minimizes $F$, we obtain for $|\varepsilon|\leq 1$
\begin{align}\label{upper_bound_F}
    F(\Psi^\varepsilon_{min})\leq F(0) = \varepsilon^2 \int_{x_l}^{x_r}q(x)^2dx\leq \int_{x_l}^{x_r}q(x)^2dx
\end{align}
Therefore, we must have $F_{-2}(\Psi^\varepsilon_{min})=0$ and $F_{-1}(\Psi^\varepsilon_{min})=0$, since otherwise $F(\Psi)$ diverges for $\varepsilon\rightarrow 0$, which would contradict Eq. (\ref{upper_bound_F}). We conclude that
\begin{align}\label{eq:delta0=0}
    \delta_0(x)=0
\end{align}
Using Eq. (\ref{eq:delta0=0}), we can restrict the set of functions to
\begin{align}
    W_h=\{w_h\in V: w_h(z,\mu)=\sum_{\nu=0}^{\infty}\varepsilon^\nu\eta_\nu(x) + \sum_{\nu=1}^{\infty}\varepsilon^\nu \delta_\nu(x),\nonumber \\ 
    ~w_h(x_l)=0~\text{for}~\mu<0,~w_h(x_r)=0~\text{for}~\mu>0\}.
\end{align}

A necessary condition for the minimum is that the derivative of $F$ with respect to $\Psi$ is zero:
\begin{align}\label{eq:F_min}
    \int_{x_l}^{x_r}\int_{-1}^{1} (L_\varepsilon\Psi^\varepsilon_{min}(x,\mu)-\varepsilon Q(x,\mu))L_\varepsilon\Psi^\varepsilon_{min}(x,\mu) d\mu dx = 0.
\end{align}
We have
\begin{align}\label{eq:L_epsilon_Psi_min}
    &L_\varepsilon\Psi^\varepsilon_{min}
    =\mu\frac{\partial\phi_0^\varepsilon}{\partial x} +\mu^2 \frac{\partial \phi_1^\varepsilon}{\partial x} +\frac{1}{\varepsilon}\mu\phi_1^\varepsilon+\varepsilon\alpha\phi_0^\varepsilon \nonumber \\
    &~~~~~~~~~~ =[\mu\eta_0'+\mu\delta_1]+\varepsilon[\mu\eta_1'+\mu^2\delta_1'+\mu\delta_2+\alpha\eta_0] +\mathcal{O}(\varepsilon^2)
\end{align}
and by inserting Eq. (\ref{eq:L_epsilon_Psi_min}) into Eq. (\ref{eq:F_min}) we obtain
\begin{align}\label{eq:F_min_epsilon}
    \int_{x_l}^{x_r}\int_{-1}^{1}\mu^2 (\eta_0' +\delta_1)^2d\mu dx +\varepsilon I_1 +\mathcal{O}(\varepsilon^2)
    = 0
\end{align}
with
\begin{align}
    &I_1=2\int_{x_l}^{x_r}\int_{-1}^{1}
    \mu[\alpha\eta_0\eta_0'+\alpha\delta_1\eta_0]+\mu^2[\eta_0'\eta_1'+\eta_0'\delta_2+\delta_1\eta_1'+\delta_1\delta_2] + \mu^3[\delta_1'\eta_0'+\delta_1\delta_1'] d\mu dx.
    \nonumber\\
    &~~~~=\int_{x_l}^{x_r}\int_{-1}^{1}2\mu^2[\eta_0'\eta_1'+\eta_0'\delta_2+\delta_1\eta_1'+\delta_1\delta_2]d\mu dx.
\end{align}
In Eq. (\ref{eq:F_min_epsilon}) $\mathcal{O}(1)$ and $\mathcal{O}(\varepsilon)$ terms show up only on the left-hand side. Therefore, both terms need to vanish, what is implied by the relation 
\begin{align}\label{eq:eta_delta_rel}
    \eta_0' = -\delta_1.
\end{align}
This shows that only terms of order $\varepsilon^2$ remain, what ends the proof.
\qed
\end{pf}

Using Lemma 1, we can prove the main theorem on PINN solutions in the diffusion limit.\\

\begin{thm}
Let the functional $F$ and the set of functions $V_h$ be as given in Eqs. (\ref{eq:min_F}) and (\ref{eq:v_h}). Suppose that $\Psi^\varepsilon_{min}$ minimizes $F$ restricted to $V_h$. Suppose further that $\varepsilon\leq 1$ and that $\Psi^\varepsilon_{min}$ has an expansion as defined in Eqs. (\ref{eq:psi_min_expansion}) and (\ref{phi_developed}).\\
Then we have: 
\begin{align}\label{eq:th1}
    \Psi^\varepsilon_{min}\rightarrow 0 \text{ as }\varepsilon\rightarrow 0~~ \text{pointwise, for all } x\in D.
\end{align}
\end{thm}

\begin{pf}
If both $\eta_\nu(x)$ and $\delta_\nu(x)$ are continuous piecewise linear functions, it follows from Eq. (\ref{eq:eta_delta_rel}) that $\eta_0$ is a linear function. With the boundary conditions we obtain $\eta_0=0$. Therefore, we have $\Psi^\varepsilon_{min}\rightarrow 0$ as $\varepsilon\rightarrow 0.$. Since all $\eta, \delta$ are continuous, the convergence holds pointwise, for all $x\in D$.\qed
\end{pf}

As a consequence, PINNs using feed forward neural networks with ReLU activation functions do not give a correct approximation of the neutron flux $\phi$ in the diffusion limit, except for $Q=0$.

\subsection{Scaling for the numerical solution}\label{subsec:scaling}

As shown in the main theorem, directly solving the neutron transport equation using a PINN with ReLU activation functions leads to a method that does not preserve the diffusion limit. Solving the neutron transport equation using first order least-squares finite elements leads to similar results (\cite{Manteuffel1998}, \cite{Varin2005}). \\
In \cite{Manteuffel1998}, a scaling of the equations that solves this issue is proposed. In the following, we also use this scaling for PINNs.\\
Recall that
\begin{align}
    P=\frac{1}{2}\int_{-1}^1d\mu,
\end{align}
Using this operator, which is a projection onto the first Legendre moment, we define the scaling operator
\begin{align}
    S=P+\tau(I-P).
\end{align}
By applying $S$ on both sides of Eq. (\ref{eq:nte_diff_limit}), we obtain
\begin{align}
    S\mu\frac{\partial \Psi}{\partial x}+\frac{\tau}{\varepsilon}(I-P)\Psi+\varepsilon\alpha P\Psi=\varepsilon Q.
\end{align}
We use the same expansion for $\Psi$ as in Eq. \ref{eq:psi_min_expansion}) to project on the first two Legendre moments. The systems of scaled equations then reads
\begin{align}\label{eq:p1_scaled}
    &\frac{1}{3}\frac{\partial \phi_1^\varepsilon}{\partial x} + \alpha\varepsilon\phi_0^\varepsilon = \varepsilon Q \nonumber\\
    & \tau \frac{\partial \phi_0^\varepsilon}{\partial x}+\tau\phi_1^\varepsilon~~ =0
\end{align}
In the following, we are only interested in the behavior of the leading orders, so we take to terms in Eqs. (\ref{phi_developed}) only up to order $\mathcal{O}(\varepsilon)$ in account. In addition, we use that $\delta_0=0$ according to Eq. (\ref{main_theorem_relations}). Therefore, we obtain
\begin{align}\label{eq:phi_leading_orders}
    &\phi_0 = \eta_0 + \varepsilon \eta_1 \nonumber\\
    &\phi_1 = \varepsilon\delta_1.
\end{align}
By inserting Eqs.~(\ref{eq:phi_leading_orders}) into Eqs.~(\ref{eq:p1_scaled}) and neglecting higher order terms, we obtain %the $P_1$-equations with the leading orders of the system:
\begin{align}\label{eq:pn_tau_scaled}
    &\varepsilon\frac{1}{3}\frac{\partial \delta_1}{\partial x} + \varepsilon\alpha\eta_0 = \varepsilon Q \nonumber\\
    & \tau \frac{\partial \eta_0}{\partial x}~~+~\tau\delta_1 =0.
\end{align}
The case $\tau=1$ corresponds to the unscaled equation, where the first equation is $\mathcal{O}(\varepsilon)$ and the second equation $\mathcal{O}(1)$, so that Eqs. (\ref{eq:pn_tau_scaled}) are unbalanced for $\epsilon\rightarrow 0$.
Choosing $\tau=\mathcal{O}(\varepsilon)$ results in a balancing of the terms in orders of $\varepsilon$. This is a common strategy to improve the convergence behavior. The numerical results in the following section demonstrate that the scaling leads to the correct diffusion limit. Here we use $\tau=\sqrt{\sigma_a/\sigma_t}=\sqrt{\alpha}\epsilon$, since this scaling directly relates to the physical parameters.

%-------------- SECTION ---------------------------------------------------------------------------------------
\section{Numerical Results}\label{sec:num-results}

In this section we present numerical results which illustrate the theory developed in Sec. \ref{sec:main-theorem}.
It is structured as follows: in Subsec. \ref{subsec:angular_discretization} we introduce the P\textsubscript{N} method which was used for the angular discretization. In Subsec. \ref{subsec:implementation} the implementation is described. In Subsec. \ref{subsec:diff_limit} and Subsec. \ref{subsec:diff_test_interface} we present numerical results for PINNs with ReLU activation function and first order LSFE to demonstrate how the scaling leads to a huge improvement in accuracy for both methods. In Subsec. \ref{subsec:tanh_test} we demonstrate that the scaling can also improve the solution for activation functions other than ReLU.

\subsection{Angular Discretization}\label{subsec:angular_discretization}

For the angular discretization of Eq. (\ref{L_diff_limit}), we use the method of moments, which is a spectral Galerkin method in $\mu$ \cite{Hesthaven_Gottlieb_Gottlieb_2007}. Let $\phi_n(x)$ be the $n$-th angular moment of $\psi$ as defined in Eq. (\ref{eq: Legendre_moments}). We then take the first $N$ Legendre polynomial moments of Eq. (\ref{eq:nte}). Using the recursion relations for Legendre polynomials, we obtain the slab geometry P\textsubscript{N} equations:
\begin{align}
    &\frac{\partial \phi_1}{\partial x}+\sigma_a\phi_0=Q\\
    &\frac{n}{2n+1}\frac{\partial \phi_{n-1}}{\partial x}+\frac{n+1}{2n+1}\frac{\partial \phi_{n+1}}{\partial x} +\sigma_n=0~(n>0).
\end{align}
We close the equations by setting
\begin{align}
    \phi_n = 0,~~~~n>N.
\end{align}
Note while we used the normed Legendre polynomials in the previous section to simplify the equations, we do not norm the Legendre polynomials here since otherwise we would not obtain the standard formulation of the P\textsubscript{N} equations in this case. In the following we use vacuum boundary conditions, which assume no incoming flux, and reflective boundary conditions, which assume that all outgoing particles are reflected back into the domain. For an $N$th order expansion where $N$ is odd, there are $(N+1)/2$ vacuum boundary conditions that read
\begin{align}
    0=\sum_{n=0}^N \frac{2n+1}{2}\phi_n(x)\int_0^{\pm 1}P_{2m-1}(\mu)P_n(\mu)d\mu~\text{for}~m=1,2,...,(N+1)/2
\end{align}
where $P_n(\mu)$ is the n-th order Legendre polynomial on the interval $[-1,1]$. Reflective boundary conditions are obtained by setting the odd moments to zero at the boundary. More details about the boundary conditions for slab geometry P\textsubscript{N} equations can be found in \cite{HamiltonEvans2015}.
It is known, that in the limit $\varepsilon\to 0$ a solution of Equation \ref{L_diff_limit} converges to the solution of a corresponding diffusion equation, which is independent of the angular variable $\mu$. Therefore, in practice a small $N$ ($N=1$ or $N=3$) is sufficient as an approximate model.

\subsection{Implementation}\label{subsec:implementation}

The PINN solver used for the numerical computations was implemented using the python library DeepXDE \cite{lu2021deepxde} with tensorflow \cite{tensorflow2015-whitepaper} as a backend. The boundary conditions were included in the loss function as an implementation of the restriction of the underlying set of functions. In order to do this, the term
\begin{align}
    \text{Loss}_{\partial D}[\Psi]=\frac{w_{\partial D}}{|T_{\partial D}|}\sum_{m=1}^{(N+1)/2}\sum_{x\in T_{\partial D}}(B_m[\Psi](x)-g(x))^2,
\end{align}
where $w_{\partial D}$ is a weight, $T_{\partial D}$ a set of training points on the boundary and $B_m$ the operator representing the $m$th boundary condition as defined in Subsec. \ref{subsec:angular_discretization}, was added to the loss function (Eq. \ref{eq:loss_function}).\\
We implemented a LSFE solver using the FEniCS library~\cite{fenics1,fenics2}. We use PETSc~\cite{petsc} as linear algebra backend in conjunction with hypre~\cite{hypre}, a library of high performance preconditioners. For this work we use a combination of a gmres solver~\cite{gmres} and an algebraic multigrid preconditioner~\cite{ruge1987algebraic}. The boundary condition was implemented by adding the term
\begin{align}
    a_b(\Psi,u):=\int_{x_l}^{x_r}\int_{\partial D} B\Psi(x,\mu)Bu(x,\mu) d\mu ds 
\end{align}
to the least-squares bilinear form (Eq. \ref{eq:LSFE_var_form}).\\
The errors in this section are computed as follows: we define $G$ points that are equidistantly distributed over the computational domain and compute the squared difference between the PINN or LSFE neutron flux and the reference solution divided by the squared reference solution. The relative error is then given by
\begin{align}
    \xi_{rel}=\sum_{g=0}^{G} \frac{(\phi_{0,PINN/LSFE}(x_g)-\phi_{0,\text{ref}}(x_g))^2}{(\phi_{0,\text{ref}}(x_g))^2}.
\end{align}
We use the zeroth flux moment since it is identical to the angular-integrated total flux. Since PINNs may converge to different solutions from different initial values \cite{Pang_2019}, we train PINNs from random initialization three times. The PINN results shown in this section are to be understood as the average of the three networks.

\subsection{Asymptotic Test of the Diffusion Limit}\label{subsec:diff_limit}

\begin{figure}[htbp]
\centering

\begin{minipage}{1.0\textwidth}
  \centering
\includegraphics[width=1.0\textwidth]{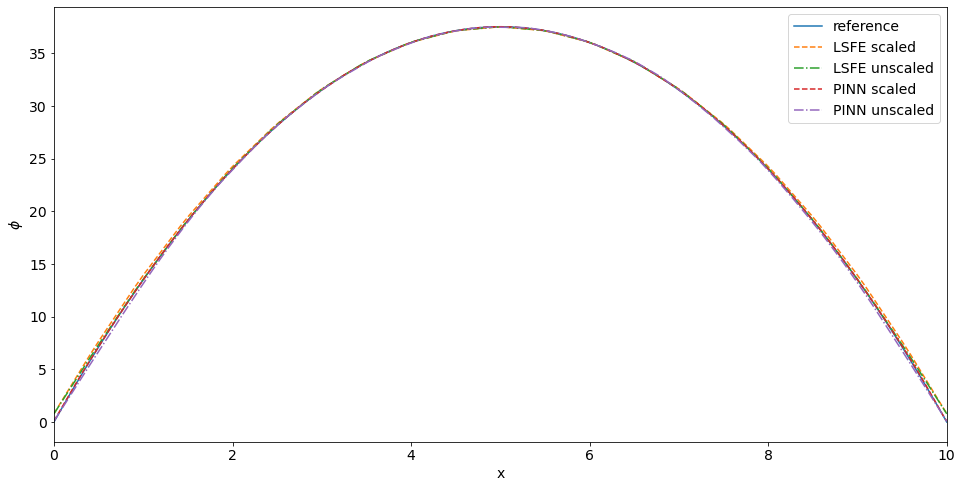}
\subcaption{$\varepsilon=10^{-2}$}\label{fig:diffusion_limit1}
\strut\end{minipage}%
\hfill\allowbreak%

\begin{minipage}{1.0\textwidth}
  \centering
\includegraphics[width=1.0\textwidth]{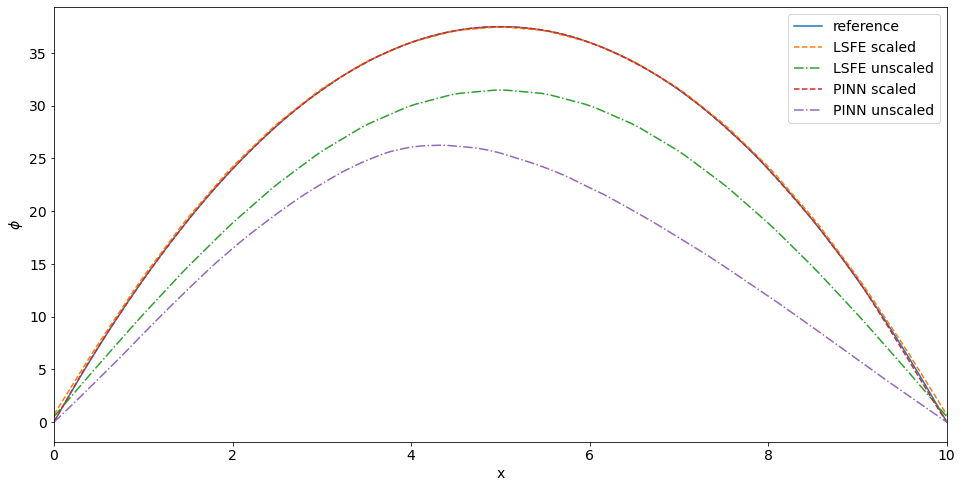}
\subcaption{$\varepsilon=10^{-3}$}\label{fig:diffusion_limit2}
\strut\end{minipage}%
\hfill\allowbreak%

\begin{minipage}{1.0\textwidth}
  \centering
\includegraphics[width=1.0\textwidth]{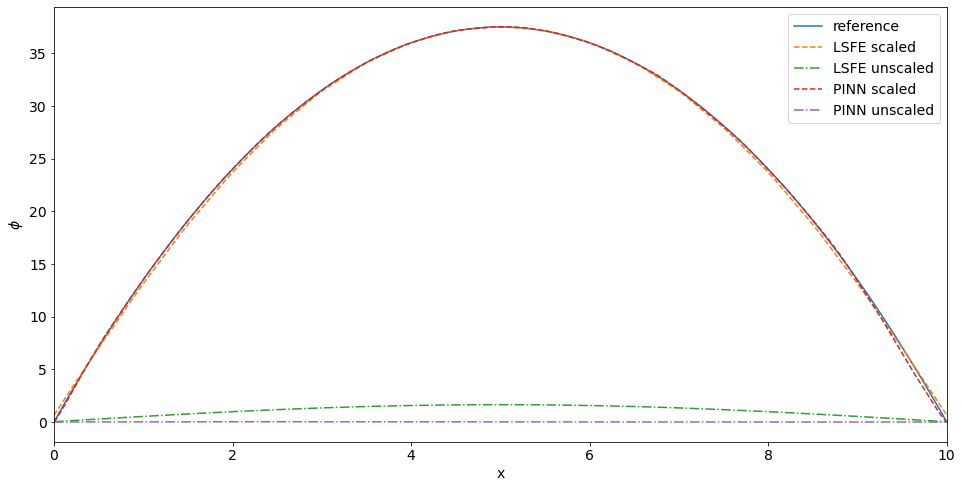}
\subcaption{$\varepsilon=10^{-4}$}\label{fig:diffusion_limit3}
\strut\end{minipage}%
\hfill\allowbreak%

\caption{Results for PINN and LSFE in an asymptotic test of the diffusion limit for three different values of $\varepsilon$.}\label{fig:diffusion_limit}

\end{figure}

The neutron transport equation in the diffusion limit is given by
\begin{align}
    \left(\mu\frac{\partial}{\partial x}+\frac{1}{\varepsilon}(I-P)+\varepsilon\alpha P\right)\Psi(x,\mu)=\varepsilon Q(x,\mu).
\end{align}
With the source term
\begin{align}
    Q(x,\mu)= 1+\alpha\phi_0(x)
\end{align}
and vacuum boundary conditions, the neutron flux $\Psi(x,\mu)$ tends asymptotically to the solution of the diffusion equation:
\begin{align}
    \phi_0(x)=-\frac{3}{2}x^2+15x.
\end{align}
In the following, we investigate the asymptotic behavior of the PINN and LSFE solutions in the diffusion limit. We choose $\alpha=10^{-2}$. It should be noted that for $\varepsilon=10^{-2}$, this choice leads to a problem described in \cite{Varin2005}.\\
For the PINN computations, we use feed forward neural networks, each with 5 hidden layers and 50 nodes per layer. An Adams optimizer with learning rate $l_r=2.5\cdot 10^{-4}$ is used for the training process. A total of 300 training points is used. For the LSFE computations, we use a mesh with 20 first order Lagrange Finite elements.\\
The results for three different values of $\varepsilon$ are depicted in Fig. \ref{fig:diffusion_limit}, while the relative errors can be found in Tab. \ref{tab:diffusion_limit}.\\
For $\varepsilon=10^{-2}$, all errors are below 1\%. The scaling reduces the already small error even further for both PINN and LSFE. For $\varepsilon=10^{-3}$, the unscaled solutions deviate significantly from the reference. While the error for the PINNs is 36.8\%, the LSFE error is 18.9\%. The scaling however reduces the error for both methods massively, the scaled PINN error is only 0.3\%, while the scaled LSFE error is 0.5\%.
For $\varepsilon=10^{-4}$, the unscaled solutions are close to zero. This is expected, since the unscaled solution is supposed to converge to zero as $\varepsilon\rightarrow 0$ according to the analysis in Sec. \ref{sec:main-theorem}. The scaled solutions however have still a small error, 0.6\% for the PINN and 1.8\% for the LSFE.\\
This results are in accordance with the theoretical predictions in Sec. \ref{sec:main-theorem}. As expected, for smaller values of $\varepsilon$ the unscaled versions of PINN and LSFE converge against zero. The scaled solution as also always better than the unscaled solution. The results for both methods are similar, however the exact numerical values deviate, as one would expect for two methods with vast differences in the implementation.

\begin{table}[!htbp]
\centering
\caption{Errors in  Fig. \ref{fig:diffusion_limit}}
\label{tab:diffusion_limit}
\begin{tabular}{|lccc|}
\hline
    $\varepsilon$ & $10^{-2}$ & $10^{-3}$ & $10^{-4}$  \\ \hline
    PINN unscaled & 0.6\% & 36.8\% & 100.0\% \\
    PINN scaled & 0.1\% & 0.3\% & 0.6\% \\
    LSFE unscaled & 0.7\% & 18.9\% & 95.8\% \\    
    LSFE scaled & 0.3\% & 0.5\% & 1.8\% \\ \hline
\end{tabular}
\end{table}

\subsection{Diffusive Test with an Interface}\label{subsec:diff_test_interface}

\begin{figure}[htbp]
\centering
\begin{minipage}{0.475\textwidth}
  \centering
\includegraphics[width=1.0\textwidth]{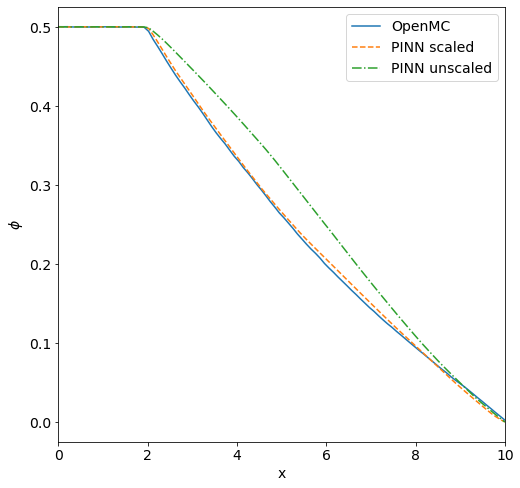}
\end{minipage}%
\begin{minipage}{0.475\textwidth}
  \centering
\includegraphics[width=1.0\textwidth]{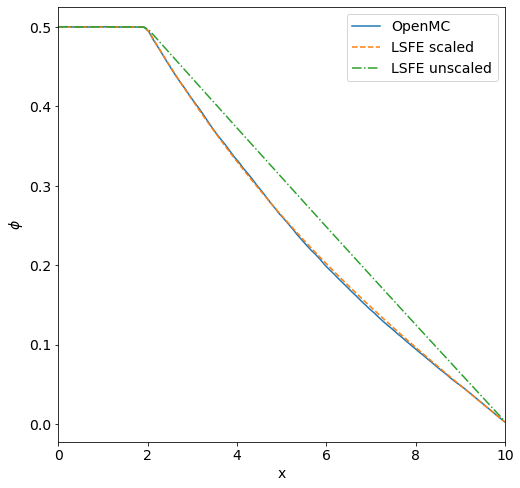}
\end{minipage}
\caption{Results for PINN and LSFE in a diffusive test with an internal interface. In both cases, the scaling leads to a huge improvement of the solution.}\label{fig:diffusive_interface}
\end{figure}

In this subsection we investigate a problem with an internal interface. The left side of the slab consists of a pure absorber $(\sigma_a=\sigma_t=2,~x\in (0,2))$, while the right side is a strong scatterer $(\sigma_t=100,~\sigma_a=10^{-4},~x\in (2,10))$ with only a very small absorption cross section. A constant source $Q=1$ is applied on the left side of the slab, i.e. for $x\in (0,2)$. On the left boundary, reflective boundary conditions are imposed, while on the right boundary vacuum boundary conditions are used. \\
For the PINN computations, we use the same architecture and learning rate as in Subsec. \ref{subsec:diff_limit}. For this testcase, 450 training points within the domain are used. For the LSFE computations we use a mesh with 50 first order Lagrange Finite elements.\\
The reference solution was computed with the OpenMC Monte Carlo Code \cite{openmc}. It should be noted that OpenMC does not approximate the neutron angle, but includes the full angular dependence, e.g. all Legendre moments. Since the higher order moments are small in the diffusive regime, we do not expect large deviations, but we still use the P\textsubscript{3}-approximation for the angular discretization to minimize potential errors.\\
Fig. \ref{fig:diffusive_interface} depicts the results for the test with an internal interface. In both cases, the scaling leads to a huge improvement of the solution. While the unscaled PINN solution deviates from the reference by 10.6\%, the scaled PINN solution deviates by only 1.4\%. The LSFE are similar. While the error for the unscaled LSFE is 9.9\%, the error of the scaled solution is only 0.6\%. We again see the similarity of the PINN und LSFE solutions.

\subsection{Diffusive Test with Hyperbolic Tangent Activation Functions}\label{subsec:tanh_test}

\begin{figure}[htbp]
\centering
  \centering
\includegraphics[width=1.0\textwidth]{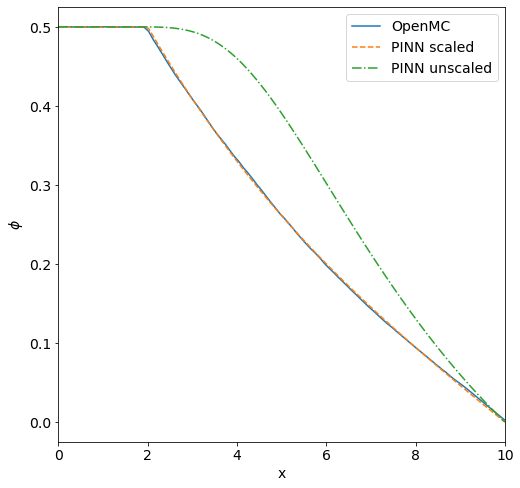}
\caption{Results for PINN with hyperbolic tangent activation functions for the diffusive test with an internal interface discussed in Subsec. \ref{subsec:diff_test_interface}. The scaling leads to a huge improvement of the solution.}\label{fig:diffusive_interface_tanh}
\end{figure}

While the analysis in Subsec. \ref{subsec:minimizer} shows that PINNs with ReLU activation functions lead to a method that does not preserve the diffusion limit, the scaling introduced in Subsec. \ref{subsec:scaling} to correct this does not assume the usage of a specific activation function. Therefore, it is close at hand to investigate the effect of the scaling on PINN with a different activation function.
In this subsection, we investigate the problem with the internal interface introduced in Subsec. \ref{subsec:diff_test_interface} for PINN with hyperbolic tangent activation functions. We use the same number of training points as before. As an optimizer L-BFGS is used since it showed better convergence properties for this testcase than the Adams optimizer. Fig. \ref{fig:diffusive_interface_tanh} depicts the results. The unscaled solution deviates form the reference by 23.4\% and is therefore worse than the unscaled solution obtained with ReLU activation functions. The error of the scaled solution is only 0.5\% and therefore smaller than in the case with ReLU activation functions. This shows that the proposed scaling leads to a vast improvement of the solution not only for ReLU activation functions, but also hyperbolic tangent.

%-------------- SECTION ---------------------------------------------------------------------------------------
\section{Discussion}\label{sec:discussion}

In this paper, we investigated the numerical stability of PINNs for multiscale transport problems. As an exemplary problem we studied the example of the neutron transport in the so-called diffusive regime. We used an analogy between PINNs and LSFE that lies in their shared approach of reformulating the solution of differential equations as the minimization of a (typically quadratic) functional. By making use of this analogy we were able to build on a theory for LFSE and adapt it for PINNs. It was shown that PINNs with ReLU activation functions yield incorrect approximate solutions in the diffusion limit. It was also demonstrated that a scaling can be applied to PINNs that leads to the correct diffusion limit. These theoretical findings were underlined with numerical results.\\

A formal proof for the convergence of the scaled solution remains a goal for future research. For LSFE, a corresponding result exists. The proof uses the fact that the quadratic functional needs to be minimized for all admissible test functions. Then, by choosing specific test functions, in addition to the relations in [Proposition 1] it can be shown, that the leading order in the formal expansion of $\phi_0^\varepsilon$ satisfies a variational form of the diffusion equation. Details can be found in \cite{Ressel_phd_thesis}. Since in the PINN context no direct analogon for the test functions exists, this cannot be mimicked with PINNs. Therefore, a different way to formally prove the convergence of the scaled PINN solution needs to be found.\\
Here, we see the possibility to build on further results of the existing theory for LSFE and adapt it for PINNs, too. For instance, there are cases known for LSFE were scalings other than the one we applied to PINNs are used. Taking additional parameters, such as the medium optical thickness or the local cell size, into account, these scalings provide a further improved accurracy. In \cite{Zheng_McLarren_2017}, as an example, such a scaling is introduced which is especially useful for optically thin materials.\\

We showed that unscaled PINNs with ReLU activation functions do not converge to the correct diffusion limit. First order LSFE exhibit the same behavior. For LSFE, the already existing theory shows that using higher order finite elements instead of a scaling can be sufficient to migitate this issue, even though applying the scaling still improves the results numerically \cite{Ressel_phd_thesis}. It is not apparent whether a PINN analogy to higher order finite elements exists and how it would look like. However, it is clear that simply choosing other activation functions is not sufficient, as our testcase with the hyperbolic tangent shows.

\section*{Acknowledgments}

The work of A. J. was performed within the project RAPID, funded by the German Federal Ministry of Education and Research (BMBF) under the grant number 033RK094B. 
The work of K.K. was performed in parts within the project RAPID, funded by the German Federal Ministry of Education and Research (BMBF) under grant number 033RK094A.
R.G.M. was supported by the International Atomic Energy Agency under contract No. USA-26472 "University of Notre Dame Contribution to the AI for Fusion Coordinated Research Project".\\
The authors are solely responsible for the content of this publication.

\pagebreak
%\bibliography{bibliography}
%\printbibliography

\bibliographystyle{elsarticle-num}
\bibliography{bibliography}
%\bibliography{\jobname}

\end{document}